\documentclass{elsart1p}
\usepackage{amssymb}

\usepackage{amsfonts}
\usepackage{graphicx}

\newtheorem{theorem}{Theorem}[section]
\newtheorem{prop.}[theorem]{Proposition}
\newtheorem{lem.}[theorem]{Lemma}
\newtheorem{cor.}[theorem]{Corollary}
\newtheorem{Bsp.}{Example}[section]
\newcommand{\norm}[2]{
\left\| #2 \right\|_{#1}
}
\newcommand{\Hil}[0]{
\mathcal{H} 
}
\newcommand{\RR}[0]{
\mathbb{R} 
}
\newcommand{\Bsp}[1]{
\begin{Bsp.} \label{#1} \bf : \end{Bsp.}
}\newcommand{\CC}[0]{
\mathbb{C} 
}
\newenvironment{proof}{\noindent \bf Proof: \rm}{\qed

\vspace{5mm}}

\begin{document}

\begin{frontmatter}
\title{Frames and Finite Dimensionality:\\ Frame Transformation, Classification and Algorithms}
\author{Peter Balazs}
\address{Austrian Academy of Sciences, Acoustic Research Institute,
         Reichsratsstrasse 17,A-1010 Vienna, Austria}
\date{\today}
\thanks{Manuscript received October ??, 20??; revised ????? ??, 2005.
        This work was partly supported by the European Union's Human Potential Programe, under contract HPRN-CT-2002-00285 (HASSIP).}
\thanks{Peter.Balazs@oeaw.ac.at}
\begin{abstract}
In this paper we will look at the connection of frames and finite dimensionality. A main focus is to present simple algorithms and make them available online. The main result is a way to 'switch' between different frames, giving an algorithm to calculate the coefficients of one frame given the analysis of another frame using their cross-Gram matrix. This is a canonical extension of the basic transformation matrix used for orthonormal bases (ONB) and is therefore called frame transformation. 
Furthermore we will summarize basic properties of frames in finite dimensional spaces. We will give basic algorithms to use with frames in finite dimensional spaces, useful for basic numerical experiments.  Finally we will give a criteria for finite dimensional spaces using frames. 
\end{abstract}
\end{frontmatter}
\maketitle


\section{Introduction} \label{intro-sec}

One mathematical background of today's signal processing algorithms, like mobile phone, UMTS, xDSL or digital television, is the concept of frames, which was introduced by Duffin and Schaefer \cite{duffschaef1}. It was made popular by Daubechies \cite{daubech1}, and today it is one of the most important foundations of Gabor \cite{feistro1}, wavelet \cite{aliant1} and sampling theory \cite{aldrgroech1}. In signal processing applications frames have received more and more attention \cite{boelc1,vettkov1}. 

To be able to work with numerically models, data and operators have to be discretized. Application and algorithms always work with finite dimensional data. There are already some works investigating frames in finite dimensional setting, for example refer to \cite{casfickinftight1,ole1}. In the finite dimensional case frames are equivalent to spanning system. Here frames are the only feasible generalization of bases, if reconstruction is wanted. In contrast to bases, frames loose the linear independence.

In this paper we will first give a short introduction to the notation of frames in Section \ref{sec:prelnot0}. In Section \ref{sec:workframfindim0} we will summarize some basic properties of frames in finite dimensional spaces. We give MATLAB \cite{matlab1} algorithms to work with frames in finite dimensional spaces. Those algorithms are either very basic and denoted here as single code lines. Or they can be downloaded from \verb_http://www.kfs.oeaw.ac.at/xxl/finiteframes/finfram1.zip_. In Section \ref{sec:frametrafo0} we will look at a way to `switch' between different frames, i.e. find a way to map between their coefficient spaces bijectively and give the algorithm for it. This is done by using the Cross-Gram matrix of the two involved frames. It is a canonical extension of the basic transformation matrix used for orthonormal bases (ONB), respectively the properties of the Gram matrix using a frame and its dual. In Section \ref{sec:claswfram0} we finish the look at the connection between frames  and finite dimensional spaces, by giving a criteria for finite dimensional spaces using frames. In particular a space is finite dimensional if and only if $\sum_k \norm{}{g_k}^2 < \infty$ for any frame $( g_k )$.


\section{Preliminaries and Notations : Frames} \label{sec:prelnot0}

 
The sequence $\left( g_k | k \in K \right)$  is called a {\em frame} \cite{duffschaef1,ole1,Casaz1,Groech1} for the (separable) Hilbert space $\Hil$, if constants $A,B > 0$ exist, such that 
\begin{equation} \label{sec:framprop1} A \cdot \norm{\Hil}{f}^2 \le \sum \limits_k \left| \left< f, g_k \right> \right|^2 \le B \cdot  \norm{\Hil}{f}^2  \ \forall \ f \in \Hil
\end{equation}
$A$ is called a {\em lower} , $B$ a {\em upper frame bound}. 
If the bounds can be chosen such that $A=B$ the frame is called {\em tight}.
%
Let $C : \Hil \rightarrow l^2 ( K )$ be the {\em analysis  operator }$ C ( f ) = \left( \left< f , g_k \right> \right) $.
Let $D :l^2( K ) \rightarrow \Hil $ be the  {\em synthesis operator}
$ D \left( \left( c_k \right) \right) = \sum \limits_k c_k \cdot g_k $. 
Let $S : \Hil  \rightarrow \Hil $ be the {\em frame  operator} 
$ S ( f  ) = \sum \limits_k  \left< f , g_k \right> \cdot g_k$. \index{frames!frame operator} \index{symbols!$S_{\mc G}$}
If it is important to note from which frame the operator was derived of, we will use an index such as for example as in $C_{g_k}$.

If we have a frame  in $\Hil$, we can find an expansion of every member of $\Hil$ with this frame.
$\left( \tilde{g}_k \right) = \left( S^{-1} g_k \right)$  is a frame with frame bounds $B^{-1}$, $A^{-1} > 0$, the so called {\em canonical dual frame}. Every $f \in \Hil$ has 
the expansions
$ f = \sum \limits_{k \in K} \left< f, \tilde{g}_k \right> g_k $
and 
$ f = \sum \limits_{k \in K} \left< f, g_k \right> \tilde{g}_k $
where both sums converge unconditionally in $\Hil$.

For two sequences in $\Hil$ $( g_k )$ and $( f_k )$, let the {\em cross-Gram matrix} $G_{g_k , f_k}$ be given by $\left( G_{g_k, f_k} \right)_{jm} = \left< f_m , g_j \right>$, $j,m \in K$.  
If $(g_k) = (f_k)$ we call this matrix the {\em Gram matrix} $G_{g_k}$. For frames $(f_k)$ the orthogonal projection $P$ from $l^2$ onto $ran(C)$ is given by
$ G ( c_k ) = \left( \left< \sum \limits_l c_l S^{-1} f_l , f_j \right> \right)_j = G_{\tilde{f}_l, f_l} c $.

\subsection{Frames in Finite Dimensional Spaces}

For a very good introduction to the topic of frames in $\CC^n$ refer to \cite{ole1} Chapter 1.

The typical properties of frames can be understood easily in the context of finite-dimensional vector spaces. It is well known \cite{ole1} that 
the finite sequences that span the whole space are exactly the finite frames. 


Let $(g_k)_{k=1}^{M}$ be a frame in $\CC^N$. The matrices 
$$ D = \left( \begin{array}{c c c c } \vline & \vline & & \vline \\
g_1 & g_2 & \dots & g_M  \\
\vline & \vline & & \vline \\
 \end{array} \right), \quad  
 C = \left( \begin{array}{c c c}\mbox{---} & \overline{g_1} & \mbox{---} \\
\mbox{---} & \overline{g_2} & \mbox{---} \\
\vdots & & \vdots \\
\mbox{---} & \overline{g_M} & \mbox{---} \\
\end{array} \right), \quad 
D^\dagger = \left( 
\begin{array}{c c c} - & \tilde{f_1} & - \\ - & \tilde{f_2} & - \\ \vdots & & \vdots \\ - & \tilde{f_k} & - \end{array}
 \right) $$
 describes the synthesis operator $D$, the analysis operator $C$ and the pseudoinverse $D^{\dagger}$. Here $D: l^2 \rightarrow V$ is defined, such that $D_{g_k} c = D \cdot \mathfrak p_{1..M}(c)$, where $\mathfrak p_{1..M}$ is the  canonical projection $l^2 \rightarrow \CC^M$. Furthermore $C: V \rightarrow \CC^M \subseteq l^2$ is the matrix such that $C_{g_k} f = C \cdot f$. And the pseudoinverse $D^{\dagger}$ of the synthesis operator is the analysis operator of the dual frame.

%
%
%
 
\section{Working with Frames in Finite Dimensional Spaces} \label{sec:workframfindim0}

Here we summarize basic properties needed later, as well as a simple algorithm to create frames.

\subsection{Examples of Infinite Frames in a Finite Dimensional Space} \label{sec:frameinffin1}

Infinite frames in finite dimensional spaces can also be constructed like in \cite{casfickinftight1}. Let us give other examples for this situation, which we refer to in Section \ref{sec:claswfram0}.

\begin{enumerate} \item Take a basis $(e_k | k = 1,..,N)$ in $\CC^N$ and let $e^{(l)}_k = \frac{1}{l} \cdot e_k$ for $l=1,2,..$. Then $(e^{(l)}_k)$ is a tight frame, as 
$$ \sum \limits_{k,l} \left| \left< f, e^{(l)}_k \right> \right|^2 = \sum \limits_{l=1}^{\infty} \sum \limits_{k=1}^{N} \left| \left< f , \frac{1}{l} \cdot e_k \right> \right|^2 = $$
$$ = \sum \limits_{l=1}^{\infty} \frac{1}{\left|l\right|^2}\sum \limits_{k=1}^{N} \left| \left< f , e_k \right> \right|^2 = \sum \limits_{l=1}^{\infty} \frac{1}{\left|l\right|^2} \norm{\Hil}{f} = \norm{\Hil}{f} \cdot \frac{\pi^2}{6}$$
\item The same is possible for $e^{(l)}_k = \frac{1}{l^2} \cdot e_k$ for $l=1,2,..$. This is again a tight frame with the bound $A= \frac{\pi^4}{90}$.
\end{enumerate}



\subsection{Frames And ONBs}

Frames can be described as images of an orthonormal basis by bounded linear operators in an infinite dimensional Hilbert space. They can even be classified by this result: 
\begin{prop.} \label{sec:framonb1}  Let $( e_k )_{k=0}^{\infty}$ be an arbitrary infinite ONB for $\Hil$. The frames for $\Hil$ are precisely the families $( U e_k )$, where $U : \Hil \rightarrow \Hil$ is a bounded and surjective operator.
\end{prop.} 

This operator is just the composition of an analysis and a synthesis operator. $U = D_{f_k} C_{e_k}$. 

This proposition seems not very easy to apply for finite dimensional applications at the first look, but it is clearly possible:


\begin{cor.} \label{sec:framonbfin1} Let $( e_k )_{k=0}^{\infty}$ be an arbitrary ONB for $l^2$. The frames for $\CC^N \subseteq l^2$ are precisely the families $( U e_k )$, where $U : l^2 \rightarrow \CC^N$ is a surjective operator.
\end{cor.} 
\begin{proof} Let $(f_k)_{k=1}^{M}$ be a frame. Define $U: l^2 \rightarrow V$ with 
$$ U(e_k) = \left\{ \begin{array}{c c} f_k & k \le M \\ 0 & otherwise \end{array} \right. .$$ 
As $(e_k)$ is an ONB this operator is well-defined and because $(f_k)$ is a frame it is surjective. It is clearly bounded. 

On the other hand $ span \left( U e_k \right) = U \left( span \left( e_k \right) \right) = U \left( l^2 \right) = \CC^N$. Therefore $(U f_k)$ is a frame.
\end{proof}

\begin{cor.} \label{sec:framinfindim1} The frames with $M$ elements in $\CC^n$ are exactly the images of an ONB 
 in $\CC^M$ by a surjective operator. 
\end{cor.} 
\begin{proof} In the above proof restrict the ONB to $\left(ker U\right)^\bot$, which in this case has dimension $M$.
\end{proof}
%
Using this corollary 
we know that matrices with full column rank correspond exactly to frames. 

\subsection{Calculation of the Frame Operator}

The frame operator 
 $S$ is a $n \times n$ matrix $S = D \cdot D^*$. This matrix can be represented very easily by
\begin{cor.} \label{sec:frammatr1} Let $( g_k )$ and $( \gamma_k )$ be sequences of elements in $\Hil$ and let $S_{g_k,\gamma_k}$ be the associated frame matrix, then
$$ {\left( S_{g_k,\gamma_k} \right)}_{m,n} = {\left( \sum \limits_{k} \gamma_k \otimes \overline{g}_k \right)}_{m,n} = \sum \limits_{k} {\left(\gamma_k \right)}_m \cdot {\left( \overline{g}_k \right)}_n $$
\end{cor.}
\begin{proof} $ S_{i,j} = \left( D_{\gamma_k} \cdot C_{g_k} \right)_{i,j} = \sum \limits_k {\left( D_{\gamma_k}\right)}_{i,k} {\left( C_{g_k}\right)}_{k,j} = \sum \limits_k (\gamma_k)_i (\overline{g}_k)_j $
\end{proof}

We can also express this term as product of an $N \times 1$ and a $1 \times N$ matrix: 
$$S_{g_k,\gamma_k} = \sum \limits_{k \in K} g_k^* \cdot \gamma_k . $$

So the frame operator matrix can be easily calculated either by multiplying $D \cdot D^*$ or by using Corollary \ref{sec:frammatr1}. 
\begin{lem.} Let $(g_k)$ and $(\gamma_k)$ be frames in $\CC^N$ with $M$ elements. Regarding numerical efficiency we see
\begin{enumerate} 
\item $S_{g_k,\gamma_k} = \sum \limits_{k \in K} g_k^* \cdot \gamma_k $ needs $M N^2 + 3 M N - 2N $ operations and
\item $S_{g_k,\gamma_k} = D_{g_k} \cdot C_{\gamma_k}$ needs $N^2 \cdot \left( 2 M - 1\right)$ operations%
\end{enumerate}
This means that the first algorithms 
is slightly more efficient (for $N \ge 4$). \\
\end{lem.}
\begin{proof}
Using the definition of matrix multiplication, it can be easily seen, that a multiplication of a $M \times N$ matrix by a $N \times P$ matrix needs $M \cdot P \cdot (2 N - 1)$ operations. This gives the result for $(ii)$. For $(i)$ we have the complex conjugation of a $N$-dimensional vector, which takes $N$ operations. The multiplication of a $N \times 1$ matrix by a $1 \times N$ matrix, which takes
$N^2$ operations. This has to be done $M$ times. Finally $M-1$ times the $N \times N$ matrices are summed, resulting in an overall sum of operations of $M \cdot \left(N + N^2\right) + \left(M - 1  \right) \cdot \left(2 N \right)$.  
\end{proof}

\subsection{Algorithmic Frame Examples} \label{algfram0}

As said before a frame can be represented as matrix with full column rank. This is just the synthesis matrix, which therefore is an easy way to represent frames. For an automatic creation of frames, we can easily create random frames by creating random matrices and checking, if the rank is full. See the file \verb_RandFrame.m_. This file uses two parameter, the dimension of the space, $dim$ and the number of frame elements $M$. In this file also for $2$- and $3$-dimensional frames the frame elements are plotted\footnote{For $3D$ plots \cite{xiong05} is used.}.

The frame bounds of a frame can be easily found by using an algorithm for the singular value decomposition, the \verb_svd_. It is well known \cite{ole1} that the frame bounds are the smallest respectively the biggest eigenvalue of the frame operator $S$. As $S = D D^*$, we get the frame bounds by using the square of the maximum and minimum of the \verb_svd_ of $D$. See for example \verb_RandFrame.m_.

\begin{figure}[!ht]
	\begin{center}
		\includegraphics[width=0.4\textwidth]{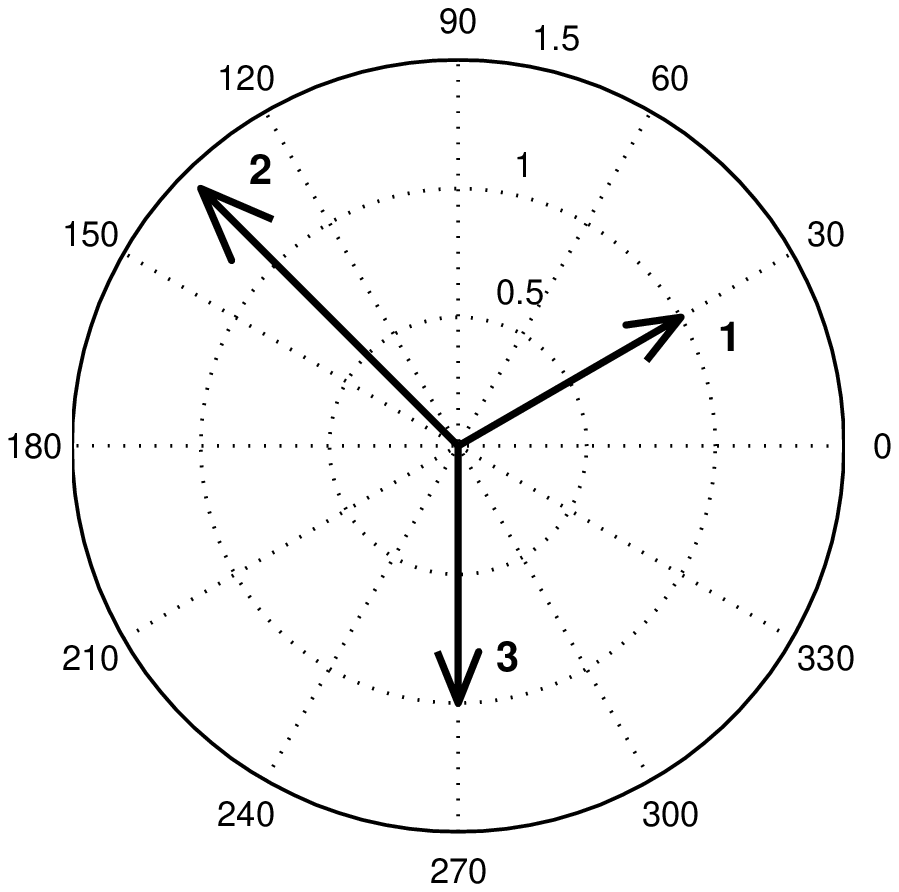} \hspace{5mm}
		\includegraphics[width=0.4\textwidth]{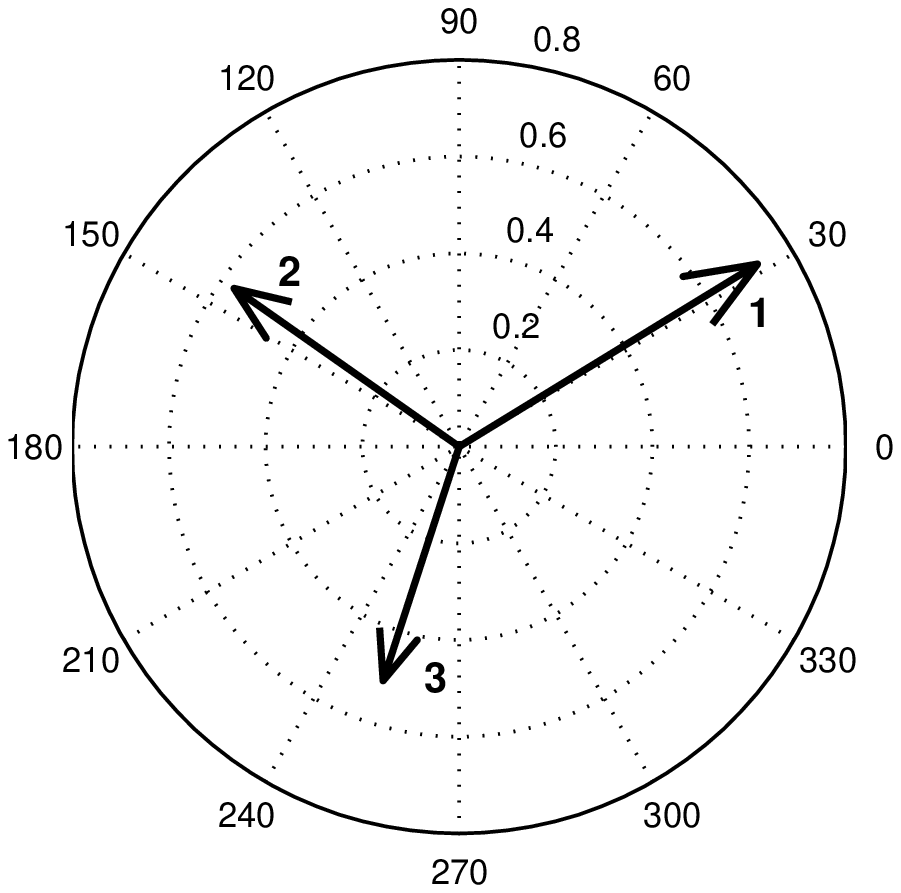}
	\end{center}
	\caption{\small \em Representation of the frame (left) and its dual (right) of Example \ref{sec:frameexamp1} $(1)$.} \label{fig:fframeexampfig1}
\end{figure}

\begin{figure}[!ht]
	\begin{center}
  	\includegraphics[width=0.4\textwidth]{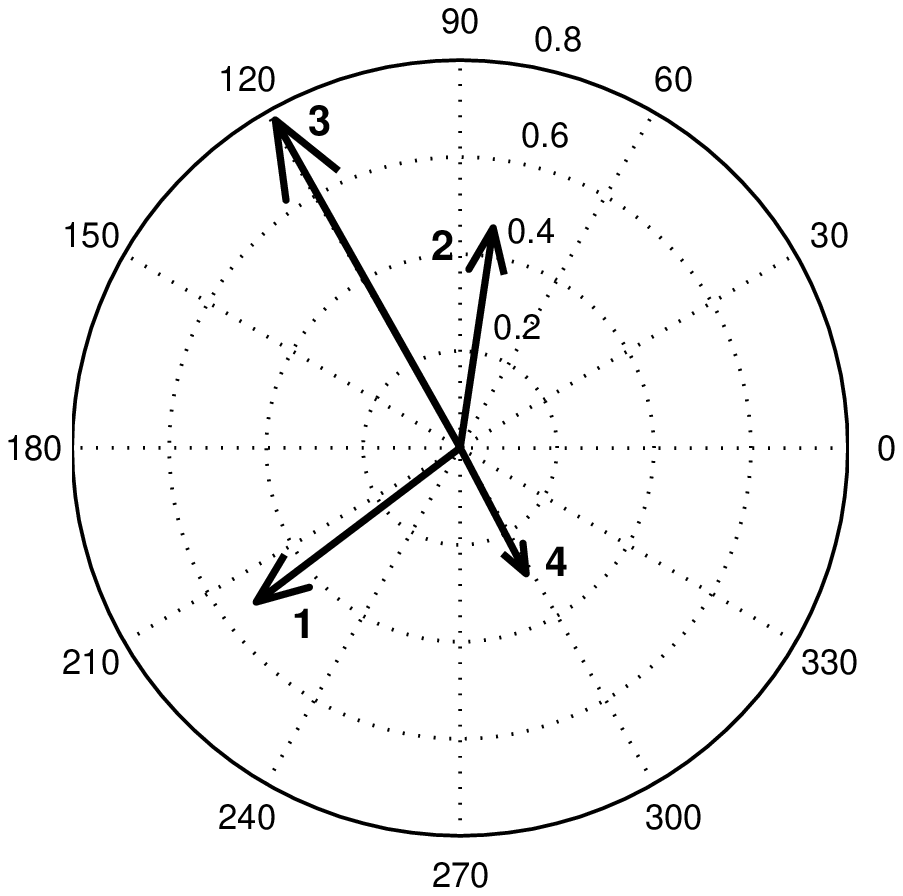} \hspace{5mm}
		\includegraphics[width=0.4\textwidth]{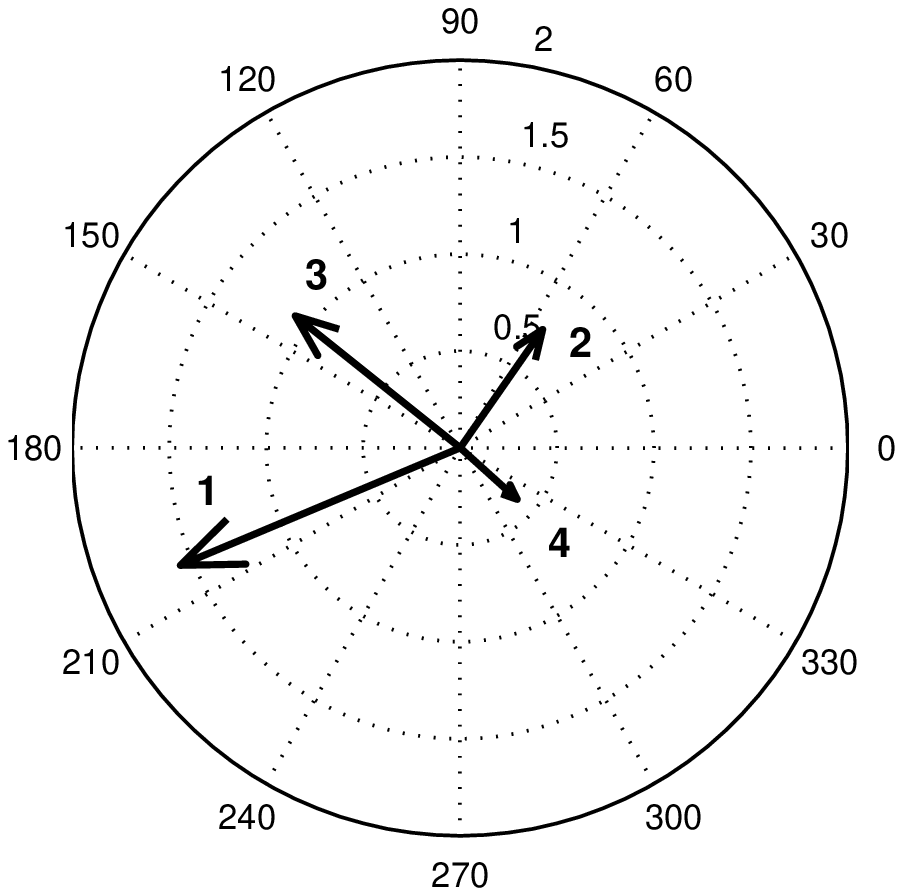}
	\end{center}
	\caption{\small \em Representation of the frame (left) and its dual (right) of Example \ref{sec:frameexamp1} $(2)$.} \label{fig:fframeexampfig2}
\end{figure}

\begin{figure}[!ht]
	\begin{center}
		\includegraphics[width=0.4\textwidth]{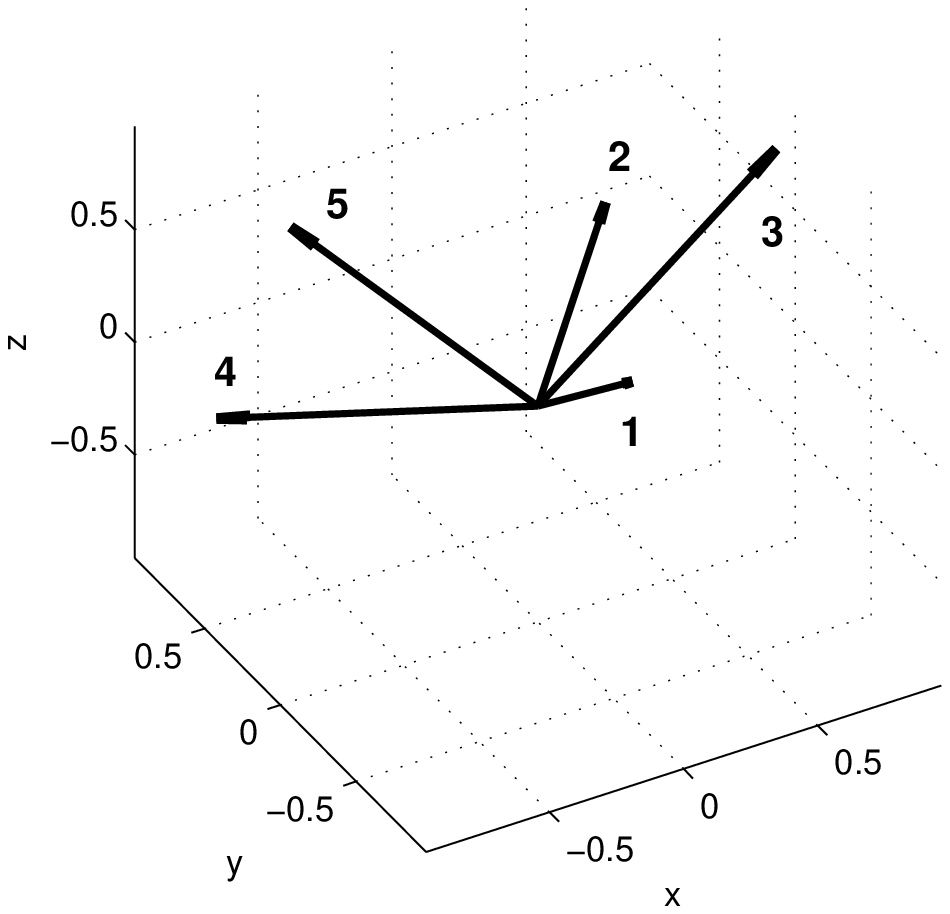} \hspace{5mm}
		\includegraphics[width=0.4\textwidth]{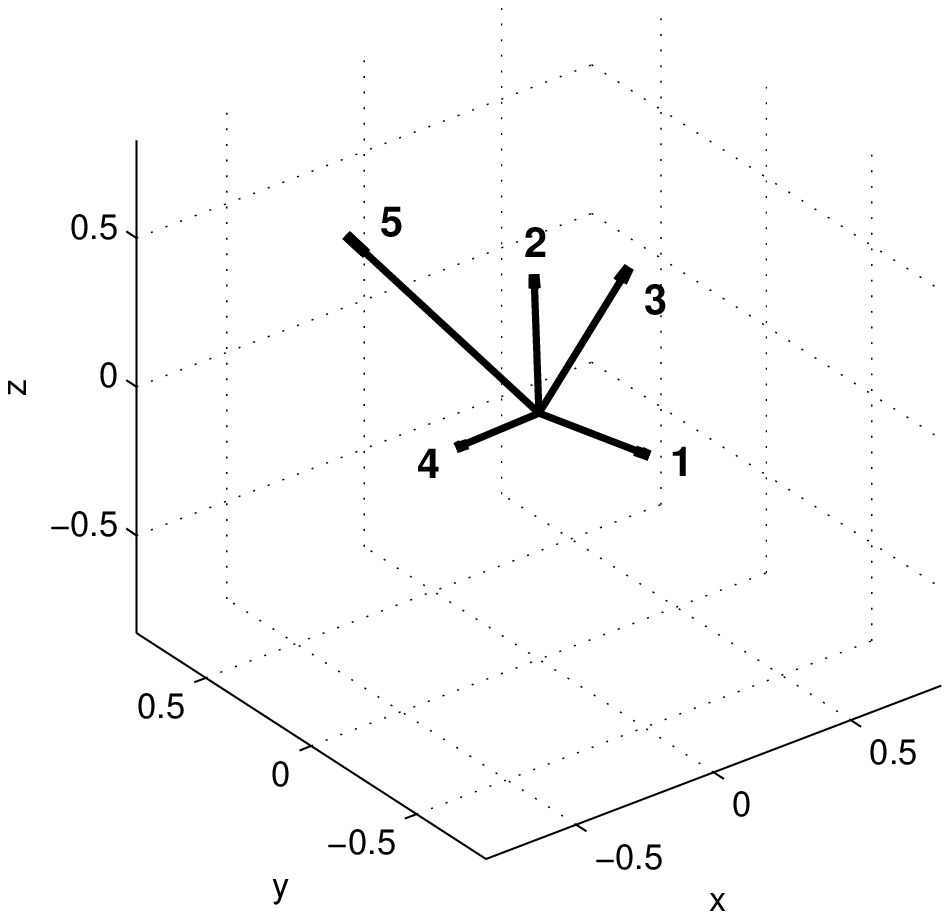}
	\end{center}
	\caption{\small \em Representation of the frame and its dual of Example \ref{sec:frameexamp1} $(3)$.} \label{fig:fframeexampfig3}
\end{figure}

\Bsp{sec:frameexamp1} 
\begin{enumerate}
\item As a non-random example let 
$$D_1 = \left( \begin{array}{r r r} \cos (30^\circ) & -1 &  0 \\ \sin (30^\circ) &  1  &-1 \end{array} \right).$$
This is a frame with bounds $A = 1.3803$, $B = 2.6197$. See Figure \ref{fig:fframeexampfig1}.
\item Using $dim = 2$ and $M = 4$ in \verb_RandFrame.m_ for example produced the matrix
$$D_2 = \left( \begin{array}{r r r r}
  -0.4205  &  0.0682 &  -0.3814  &  0.1361 \\
   -0.3176  &  0.4542 &   0.6770  & -0.2592 
 \end{array} \right) .$$
 This frame has the frame bounds $A = 0.313485$ and $B = 0.864726$. See Figure \ref{fig:fframeexampfig2}.
\item  For $dim = 3$ and $M = 5$ we got as an example  
$$D_3 = \left( \begin{array}{r r r r r} 
-0.9803 &  -0.6026 & -0.6024 &  -0.1098 &  -0.1627 \\
   -0.7222  &  0.2076 &  -0.9695 &   0.8636  &  0.6924 \\
   -0.5945 &  -0.4556 &   0.4936 &  -0.0680  &  0.0503
 \end{array} \right) .$$
       as an example of a frame with $5$ elements in $\RR^3$. It has the frame bounds $A = 0.4660$ and $B = 2.4484$. (For plotting of $3$-dimensional frames see \verb_Plot3DFrame.m_.) See Figure \ref{fig:fframeexampfig3}.
\end{enumerate}

The dual frame can be calculated using an algorithm for the pseudoinverse and typing \verb_Cd = pinv(D)_, where $D$ is the synthesis matrix of the original frame. 

\section{Frame Transformation} \label{sec:frametrafo0}

Linear algebra tells us that invertible matrices are exactly describing the transformation from one bases to another and unitary operators are exactly representing the change between orthonormal bases. So "switching" from one ONB to another is rather straight forward. But what about frames? How can a representation be changed from one frame to another in the finite dimensional case? 

Let $(g_k)_{k=1}^M$ and $(f_i)_{i=1}^N$ be two frames. We want to find a way to switch between the coefficients of these frames. 
We extend a well-known property of the Gram matrix of a frame and its dual to arrive at
\begin{theorem} \label{sec:framtrans1} Let $(g_k)_{k=1}^M$ and $(f_i)_{i=1}^N$ be two frames. The $M \times N$ matrix $G = G_{f_j,\tilde{g}_k}$ maps $ran ( C_{g_k} )$ bijectively onto $ran ( C_{f_k} )$ such that 
$$ \left( \left<f,g_k\right> \right)_k \mapsto \left( \left< f , f_k \right> \right)_k $$
and
$$ f = \sum \limits_{i=1}^M \left<f,g_i\right> \tilde{g}_i = \sum \limits_{i=1}^N \left( G \cdot \left( \left< f , g_k \right> \right)_{k=1}^M \right)_i \tilde{f}_i $$
\end{theorem}
\begin{proof} \footnote{This could also be shown using properties of the matrix description of operators using frames \cite{xxlframoper1}.}
Let $c \in ran(C_{g_k})$, there exists $f \in \Hil_1$ such that $c_k = \left< f, g_k \right>$. Then 
$$ \left( G \cdot c \right)_i = \sum \limits_k G_{i,k} c_k = \sum \limits_k \left< \tilde{g}_k , f_i \right> \cdot  \left< f, g_k \right>  = \left< \sum \limits_k \left< f, g_k \right> \tilde{g}_k , f_i \right> = \left< f , f_i \right> .$$  
Given $c \in ran ( C_{f_k} )$ for $f = \sum \limits c_k \tilde{f}_k$ the element $d = \left< f , g_k \right>$ is mapped on $c$, as
$$ \left( G \cdot d \right)_n = \sum \limits_m G_{n,m} d_m = \sum \limits_m \left< \tilde{g}_m , f_n \right> \cdot \left< f  , g_m \right> = $$
$$ = \left< f , f_n \right> = \sum \limits_k c_k \left< \tilde{f}_k, f_n \right> = G_{f_n,\tilde{f}_k} \cdot c = c$$
because $c \in ran ( C_{f_k} )$ and $G_{f_n,\tilde{f}_k}$ is the projection on this space. So this mapping is surjective.

Now suppose $G \cdot c = G \cdot d$ for $c = \left< f , g_k \right>$ and $d = \left< g , g_k \right>$. This means that 
$$ \left< f , f_i \right> = \left< g , f_i \right> \Longleftrightarrow f = g \Longleftrightarrow c = d $$
and so the mapping is injective and well-defined. 
\end{proof}

So analogue to the basis transformation matrix defined in linear algebra, the Gram matrix $G = G_{f_j,\tilde{g}_k}$can be called the {\em frame transformation matrix}

\subsection{Examples and Algorithms}

We apply this result on the examples given above. We want to do the frame transformation from the frame in Example \ref{sec:frameexamp1} $(1)$ to the one in Example \ref{sec:frameexamp1} $(2)$: 

The cross Gram matrix $G = G_{f_j,\tilde{g}_k}$ can be calculated very easily. 
As a matrix clearly 
$$G_{m ,j } = \left< \tilde{g}_j , f_m \right> = \left( D_{f_j}^* \cdot D_{\tilde{g}_k} \right)_{m,j} = D_{f_j}^* \cdot \left( D_{g_k}^* \right)^\dagger . $$

In MATLAB the above (easy) calculation of the Cross-Gram matrix is done in the file \verb_CrossGram.m_: 
\begin{center}
\verb_G12 = CrossGram(pinv(D1)',D2) = (D2')*(pinv(D1)');_ 
\end{center}

This can be applied on a random vector $f$ to show that the coefficients $ \left<f,g_i\right>$ are mapped to $\left< f , f_k \right>$. See the file \verb_testframtrans.m_ . To test the numerical stability of this transformation, this is repeated $100'000$ times. The maximum error, i.e. the norm of the difference of the coefficients with the second frame minus the Gram matrix times the coefficients of the first frame, is around $100 \cdot \epsilon$, where $\epsilon$ is the smallest non-zero error possible using the given precision.
\\

This can also be used together with the algorithm \verb_RandFrame.m_ to compare two randomly created frames. Also in this case it can be nicely seen using the methods above, that the error is around the precision error. See \verb_testframtrans2.m_ . An exemplary output of this algorithm would be:

{\footnotesize
\begin{verbatim}
Dimension = 8, Number of frame elements, first frame = 19, second frame = 27
First frame:
Created a frame with lower bound A = 2.37804 and upper bound B = 11.149
Second frame:
Created a frame with lower bound A = 2.75169 and upper bound B = 20.3146
The maximal error in 100 runs was 2.01704e-014
\end{verbatim}
}

\section{Classification Of Finite Dimensional Spaces With Frames} \label{sec:claswfram0}

For an ONB $(e_i)$ if the sum of the elements $\sum \limits_i \norm{\Hil}{e_i}$ is finite, the dimension of the space is finite and vice versa\footnote{This is equivalent to the identity being compact. Also this can be extended to the frame case, as it is easy to show: The frame operator $S$ is compact if and only if $\Hil$ is finite dimensional.}. The Example \ref{sec:frameinffin1} shows that that is not true anymore with frames, as in this case 
$$ \sum \limits_{l,k} \norm{}{e^{(l)}_k} = \sum_{l=1}^\infty \sum_{k=1}^N \norm{}{\frac{1}{l} e_k} = \sum_{l=1}^\infty \frac{N}{|l|} = \infty$$ 

But taking the square sum of the norms of the elements of a frame for $\Hil$ is an equivalent condition for $\Hil$ being finite dimensional:
\begin{prop.} \label{sec:frafindim1} Let $(g_k)$ be a frame for the Hilbert space $\Hil$. Let $(e_l)$ be an ONB for $\Hil$. Then the following statements are equivalent
\begin{itemize}
\item $\sum_k \norm{}{g_k}^2 < \infty$
\item $\sum_l \norm{}{e_l}^2 < \infty$
\item the space is finite dimensional.
\end{itemize} \end{prop.} 
\begin{proof} The equivalence of the second and third statements is clear. 
$$\sum_k \norm{}{f_k}^2 = \sum_k \sum_l \left| \left< f_k, e_l\right> \right|^2 
 $$
On the one hand, when first sum is finite, 
the summation can be interchanged and this means
$$\sum_k \norm{}{f_k}^2 = \ge \sum_l \sum_k \left| \left< f_k, e_l\right> \right|^2  \ge \sum_l A \cdot \norm{}{e_l}^2 $$
and so the sum $\sum_k \norm{}{e_k}^2$ must be finite.

On the other hand when $\sum_l \norm{}{e_l}^2 < \infty$ then $\sum_l \sum_k \left| \left< f_k, e_l\right> \right|^2  \le \sum_l B \cdot \norm{}{e_l}$ and so the sum is convergent and the order of the summation can be changed and therefore
$$\sum_k \norm{}{f_k}^2  \le \sum_l B \cdot \norm{}{e_l} $$
\end{proof}

And as an evident corollaries we find:
\begin{cor.} \label{sec:framdim1} Let $(f_k)$ be a frame and $(e_k)$ an ONB for $\Hil$ then
$$ A \cdot \sum_l \norm{}{e_l}^2 \le \sum_k \norm{}{f_k}^2 \le B \cdot \sum_l \norm{}{e_l}^2 $$
or equivalently (for finite dimensional spaces)
$$ A \cdot dim \left( \Hil \right) \le \sum_k \norm{}{f_k}^2 \le B \cdot dim \left( \Hil \right) $$
\end{cor.}

\begin{cor.} Let $(f_k)_{k=1}^m$ be a tight frame in the finite dimensional $\Hil$ with $dim \Hil = n$, then
$$ \sum_k \norm{}{f_k}^2 = A \cdot n \quad \mbox{ resp. } \quad \frac{\sum_k \norm{}{f_k}^2}{n} = A $$
If all frame elements have equal length, i.e.  
$\norm{\Hil}{f_k} = d$ for all $k$, then
$$ m \cdot d = A \cdot n \quad \mbox{ resp. } \quad \frac{m \cdot d }{n} = A $$
If this frame is normalized, i.e. $d = 1$, then
$$ A = \frac{m}{n}$$
\end{cor.}
Compare to \cite{reawal02}, where a possibility to construct such a frame is given.


%
%

\section{Algorithms}

All algorithms can be found at \\ \verb_http://www.kfs.oeaw.ac.at/xxl/finiteframes/finfram1.zip_

This includes the files:
{\small 
\begin{itemize}
\item \verb_CrossGram.m_
\item \verb_Plot3DFrame.m_
\item \verb_RandFrame.m_
\item \verb_testframtrans.m_
\item \verb_testframtrans2.m_
\end{itemize}
}

\small

\end{document}